\documentclass[12pt]{article} 
\usepackage{amsmath}
\usepackage[dvips]{graphicx}
\usepackage{amsfonts}
\usepackage{amsopn}
\usepackage{amsthm}
\def\bb{\bf}

\def\spt{{\it spt}}
\def\d{d+1}

\def\Reals{{\bb R}}

\def\RK{{\Reals}^{K}}

\def\Rn{{\Reals}^{n+1}}
\def\Rnn{{\Reals}^{n}}

\def\Gam{\Gamma(\mu,\s)}
\def\Gamp{\Gamma_+(\mu,\s)}

\def\xN{{\langle x, N \rangle}}
\def\xiN{{\langle x_i, N \rangle}}

\def\calA{{\cal A}}

\def\calF{{\mathcal F}^n}
\def\calFK{{\calF_K}}
\def\bFK{{\bar{\cal {F}}}_K^n}

\def\calQ{{\cal Q}}
\def\calQK{{\cal Q}_K}

\def\calC{{\cal C}}

\def\hK{{\hat {K}}}

\def\hK{{\hat{\cal K}}}

\def\calO{{\cal O}}

\newcommand{\Sn}{{\bf S}^n}
\newcommand\Snn{ {\bf S}^{n-1}}

\newcommand{\ep}{\epsilon}
 
\newcommand{\p}{\rho}

\newcommand{\n}{\nabla}

\newcommand{\Ht}{\tilde {h}}
\newcommand{\Pt}{\tilde {P}}
\newcommand{\Ft}{\tilde {F}}
\newcommand{\Pd}{ P^{\d}}
\newcommand{\pt}{\tilde {\p}}
\newcommand{\hpt}{\hat{\pt}}
\newcommand{\hHt}{\hat{\Ht}}

\newcommand{\pd}{{\p}^{\delta}}
\newcommand{\hpd}{{\hat{\p}}^{\delta}}
\newcommand{\hd}{h^{\delta}}
\newcommand{\hhd}{{\hat{h}}^{\delta}}

\newcommand{\ox}{\bar {x}}

\newcommand{\hp}{\hat {\p}}
\newcommand{\hh}{\hat {h}}

\renewcommand{\o}{\omega}
\newcommand{\s}{\sigma}
\newcommand{\g}{\gamma}
\newcommand{\gt}{\tilde{\gamma}}
\renewcommand{\a}{\alpha}

\renewcommand{\d}{\delta}

\newcommand{\ra}{\rangle}

\def\th{\theta}
\def\Th{\Theta}

\def\calA{{\cal A}}
\def\calC{{\cal C}}

\def\tg{{\tilde{\gamma}}}

\newcommand{\la}{\langle}

\newtheorem{theorem}{Theorem}
\newtheorem{definition}[theorem]{Definition}

\newtheorem{corollary}[theorem]{Corollary}

\newtheorem{proposition}[theorem]{Proposition}
\newtheorem{remark}[theorem]{Remark}
\title {Embedding $\Sn$ into $\Rn$ with
given integral Gauss curvature and 
 optimal mass transport on $\Sn$} 
\author{Vladimir Oliker\thanks{The research of the author
was partially supported by the National Science Foundation grant DMS-04-05622
and the Air Force Office of Scientific Research (AFOSR) under contracts FA9550-05-C-0058 and  FA9550-06-C-0058.
}\\
Department of Mathematics and Computer Science,\\
 Emory University, Atlanta, Georgia\\
oliker@mathcs.emory.edu}
\date{}
\begin{document}
\maketitle
\pagenumbering{arabic}
\setcounter{section}{0}
\setcounter{subsection}{0}
\begin{abstract}  In \cite{Al_pol}
A.D. Aleksandrov raised a general question of finding variational statements and proofs of existence of polytopes with given geometric data. The first goal of this paper is to give a variational solution to the problem of existence and uniqueness of a closed convex
hypersurface in Euclidean space with prescribed integral Gauss curvature. Our solution includes the case of a convex polytope. This
problem was also first considered by Aleksandrov and below it is referred to as Aleksandrov's problem. The second goal of this paper is to show that in variational form the Aleksandrov problem
is closely connected with the theory of optimal mass transport on a sphere with  cost function and constraints arising naturally from geometric considerations.

\end{abstract}
\section{Introduction}\label{intro}
In his book on convex polytopes A.D. Aleksandrov raised a general
question of finding variational statements and proofs of existence of polytopes 
with given geometric data \cite{Al_pol}, section 7.2.4. As examples of geometric problems in convexity theory for which variational solutions are possible Aleksandrov referred to the Minkowski problem for which such a proof was given by Minkowski himself \cite{Minkowski:03}  (see, also, Schneider \cite{Schneider}, section 7.1 and Klain \cite{Klain:2004}) and the Weyl problem for which a variational approach was sketched by Blaschke and Herglotz \cite{Blaschke-Herglotz} and later carried out by Volkov \cite{volkov:60-2}.

The first goal of this paper is to provide a variational solution to the problem of existence and uniqueness of a closed convex
hypersurface in Euclidean space with prescribed integral Gauss curvature. Our solution includes the case of a convex polytope.
This problem was first considered also by Aleksandrov  and we refer to it as Aleksandrov's problem. Aleksandrov studied it with a
nonvariational approach in \cite{Al_invariance} (see also \cite{Al_pol},
ch. 9) using his so called ``mapping
lemma'' which is a variant of the domain invariance theorem and  requires uniqueness of the solution in order to establish existence.
By contrast, our solution of Aleksandrov's problem, in addition to being variational, does not require uniqueness for
the proof of existence. 

The second goal of this paper is to show that in variational form Aleksandrov's problem 
can be considered as a problem of optimal mass transport on a sphere
with appropriate cost function and constraints.  As a result of this connection with optimal transport ``economics-like'' interpretations of some classical geometric concepts such as the Gauss map are obtained. Furthermore, this connection permits,
in principle, a 
numerical determination of a convex polytope by the integral
Gauss curvature by methods of linear programming. This should be useful in
some applied problems involving numerical determination of convex surfaces by Gauss curvature \cite{Little, La:1993, winston, oliker:spie:5942}. 

In order to state the Aleksandrov problem and outline our main results we recall the notion of integral Gauss curvature.
In Euclidean space $\Rn,~n\geq 1,$ fix a Cartesian coordinate system
with the origin at some point $\calO$. By a closed convex hypersurface
we understand here the boundary of a compact convex set in $\Rn$ containing $\calO$ in its interior. The set of such closed convex hypersurfaces is denoted by $\calF$. 
 Let 
$\Sn$ be a sphere of unit radius centered at $\calO$ and $F\in\calF$.
Since $F$ is star-shaped with respect to $\calO$, we can parametrize it as $r(x)=\p(x)x, ~x \in \Sn$, where
$\p(x)$ is the radial function giving the distance from $\calO$
to the point of intersection of $F$ with the ray of direction $x$ originating at 
$\calO$. (Here and elsewhere in the paper a point on $\Sn$ is treated also as a unit vector in $\Rn$ originating at $\calO$.)  The Gauss map $N_F:F\rightarrow \Sn$
maps a point $X\in F$ into $N_F(X)\subset \Sn$, which is the set of outward unit normals to all supporting hyperplanes to
$F$ at $X$. Define the {\it generalized Gauss map} as $\a_F=N_F\circ r:\Sn \rightarrow \Sn$. At $x\in\Sn$ such that at  $r(x)\in F$ there is more than
one supporting hyperplane the map $\a_F$ is multivalued.  If $\o$ is a subset  
of $\Sn$ then  we set 
$\a_F(\o)=\bigcup_{x\in \o}\a_F(x)$. For any Borel subset $\o\subset\Sn$ the set $\a_F(\o)$ is Lebesgue measurable on $\Sn$
\cite{Schneider}, section 2.2. The function $\s(\a_F)$,
where $\s$ is the standard $n-$dimensional Lebesgue measure on $\Sn$,  is a finite, nonnegative and countably additive measure on Borel subsets of $\Sn$.
 It
is called the integral Gauss curvature of $F$ (transferred to $\Sn$) \cite{Al_invariance}.

The Aleksandrov problem posed by him in \cite{Al_invariance}  is to
find conditions on a given measure $\mu$ on $\Sn$ under which there exists an $F \in \calF$ such that 
\begin{equation}\label{A-eq}
\s(\a_F(\o)) = \mu(\o) ~\mbox{ for any Borel set}~\o \subset\Sn. 
\end{equation}
In the same paper Aleksandrov proved the following
\begin{theorem}\label{AD1} In order for a given function $\mu$ on Borel subsets of
$\Sn$ to be the integral Gauss curvature of a convex hypersurface $F \in \calF$
it is necessary and sufficient that\\
\indent{\em (i)} $\mu$ is nonnegative and countably additive on Borel 
subsets of $\Sn$;\\
\indent {\em (ii)} $\mu(\Sn) = \s(\Sn)$;\\
\indent {\em (iii)} the inequality  
\begin{equation} \label{eq1}
\mu(\Sn \setminus \o) > \s(\o^*),
\end{equation}
holds for any spherically convex subset $\o \subset \Sn, ~\o\neq \Sn$; here, $
\o^*$ is the dual to $\o$, that is, $\o^* = \{y \in \Sn~|~ \la x, y 
\ra \leq 0 ~\forall x \in \o\}$.
\end{theorem}
The necessity of (i)-(iii) is easy to verify. The proof of sufficiency given by Aleksandrov consists of two steps. In step
one the theorem is established for convex polytopes when the function $\mu$ is an atomic measure concentrated
at a finite number of points $x_1,...,x_K \in \Sn$. Here,
to prove existence Aleksandrov uses his mapping lemma proving first uniqueness among polytopes  up
to a homothety with respect to $\calO$.
In step two the initially given function $\mu$ is approximated weakly by a sequence
of atomic measures for which, by step one, there exists a sequence of polytopes in $\calF$ solving the problem in this case. It is shown
that this sequence contains a subsequence converging to 
a convex hypersurface in the Hausdorff metric, the limiting hypersurface is in $\calF$ and its integral Gauss curvature  is  $\mu$. 
 Uniqueness in $\calF$ up
to a homothety with respect to $\calO$ was proved by Aleksandrov in  \cite{Al:42}. Later, Pogorelov
\cite{P}, \cite{P1}, \cite{P3} studied the Aleksandrov problem
and its generalizations by a different approach not relying on the
mapping lemma. The variational solution presented below is 
different from that by Pogorelov and gives significantly more information on the solution and its properties.

In this paper the hypersurface required in Theorem \ref{AD1}  is obtained as a minimizer for 
the functional
\begin{equation}\label{dual0}
\calQ[h,\p] = \int_{\Sn}\log h(N)d\s(N) - \int_{\Sn}\log \p(x)d\mu(x)
\end{equation}
on the set
\begin{eqnarray}\label{dual1}
\calA=\{(h, \p) \in C(\Sn)\times C(\Sn)~| ~h > 0,~\p > 0,\\ \nonumber
~\log h(N)-\log \p(x) \geq \log \xN,~(x,N) \in \Sn\times \Sn,~\xN > 0\}.
\end{eqnarray}
It is shown that a minimizing pair $(\Ht,\pt)$ exists, unique (up to
a multiplicative constant) and $\Ht$ and $\pt$ are,
respectively,
the support and radial functions of a hypersurface $\Ft\in \calF$ satisfying (\ref{A-eq}). 
In agreement with
the question of Aleksandrov asking for a variational solution in the class of convex polytopes we first solve the minimization problem
for polytopes and then prove
existence of minimizers in the general case. As can be seen from considerations in section \ref{var1}, the choice of the  functional $\calQ$ and
of the set $\calA$ is a completely natural consequence of the geometric duality relating the radial and support functions of a convex hypersurface in 
$\calF$ via a Legendre-like transform.   

It turns out that the generalized Gauss map $\a_{\Ft}$ of the minimizer of $\calQ$ maximizes the functional
\begin{equation}\label{monge0}
\int_{\Sn} c(\th^{-1}(N),N)d\s(N),~\th \in \Th,
\end{equation}\
where the ``cost'' function is given for $(x,N) \in \Sn\times \Sn$ by
\begin{equation}\label{cost}
c(x,N) = \left \{ \begin{array}{ll}
\log \xN  & \mbox{when~ $\xN > 0$}\\
-\infty & \mbox{otherwise}
\end{array}
\right .
\end{equation}
and $\Theta=\{\theta\}$ is the set of measurable maps of $\Sn$ onto itself, possibly multivalued,  satisfying the following conditions:
\begin{eqnarray}
\mbox{the image $\theta(\o)$ of any Borel set $\o \subset \Sn$ is
Lebesgue measurable,} \label{theta1} \\
\s\{N\in \Sn ~|~ \theta^{-1}(N) ~\mbox{contains more than one point} \}=0,
\label{theta2}
\end{eqnarray}
and for $\mu$ as in Theorem \ref{AD1} and any $f \in C(\Sn)$
\begin{equation}\label{g4}
\int_{\Sn}f(\th^{-1}(N)) d \s(N)  = \int_{\Sn}f(x)d\mu(x).
\end{equation}
In addition, the maximizer $\a_{\Ft}$ is unique (up to sets of measure zero with respect to $\mu$).

From the geometric point of view, the function $\log \la x, N_F(x)\ra$ gives a scale invariant quantitative measure of ``asphericity'' of a hypersurface $F$ with respect to $\calO$. For example, for a sphere centered at $\calO$  it is identically zero, while for a sufficiently elongated ellipsoid of revolution centered at $\calO$ it has large negative values at points where the radial direction is nearly orthogonal to the normal. Roughly speaking, the above result says that the most efficient way (with respect to the cost $c(x,N)$) to transfer to $\s$ an abstractly given measure $\mu$ on $\Sn$ is to move it by the least possible distance and this is accomplished by the generalized Gauss map $\a_{\Ft}$ of the convex hypersurface $\Ft$ solving the problem (\ref{dual0}), (\ref{dual1}).

The above result establishing that $\a_{\Ft}$ maximizes the functional in (\ref{monge0}) can also be viewed
as a counterpart of a result by Brenier \cite{Brenier}, \cite{Brenier1} (and for more general cost functions by Gangbo-McCann \cite{Gangbo/McCann:95}, \cite{Gangbo/McCann:96} and
Caffarelli \cite{Caf_alloc:96})
showing that among measure preserving maps between two 
sets $U$ and $V$ in Euclidean space the map optimal with respect to the cost function $|x-y|^2,~x \in U,~y\in V$, is unique and is the subgradient of a convex function.


In the framework of mass transport theory on $\Sn$ the problem  (\ref{dual0}), (\ref{dual1}) is the dual of the following primal 
maximization problem of finding $\gt$ such that
\begin{equation}\label{primal}
\calC[\gt]=\sup_{\Gam}\calC[\g],~\mbox{where}~\calC[\g]:=\int_{\Sn}\int_{\Sn}c(x,N) d \g(x,N),
\end{equation}
and $\Gam$ is a set of joint Borel measures on $\Sn \times \Sn$ with
marginals $\mu$ and $\s$ as in Theorem \ref{AD1}. In this case, the optimal measure is generated by
the map $\a_{\Ft}$, that is, $\gt[U,V] =\s[\a_{\Ft}(U)\cap V]$ for any Borel subsets $U$ and $V$ on 
$\Sn$. In addition, the duality relation  
\[
\calQ[\Ht,\pt]=\calC[\gt]
\]
 holds.
Thus, similar to the classical theory of optimal mass transport developed
by L.V. Kantorovich \cite{KA:39}, \cite{KA} and many other
authors, we also have here the primal  and dual 
 problems  and the usual duality between them. In this respect, the problem (\ref{primal}) can be viewed as a direct generalization of the classical problem. On the other hand, the problem (\ref{primal}) is essentially different from the classical problem. In our case the transport takes
place on a sphere and the cost function (\ref{cost}) does not have a
direct analogue among cost functions considered previously in Euclidean space \cite{Gangbo/McCann:95}, \cite{Gangbo/McCann:96}, \cite{Caf_alloc:96} or in
Riemannian space \cite{feldman-mccann-tams:02}. 
Cost functions of the type (\ref{cost}) seem to be more common (and natural) in geometric problems
than in economics or mechanics.

It is worthwhile noting that the variational approach used in this paper  is quite general and is expected to be useful for studying other nonlinear geometric problems, in particular,
the generalized Minkowski problem considered in \cite{LutOli}, \cite{Lut-Yang-Zhang:04},
\cite{Lut-Yang-Zhang:06}.
In a simpler situation this approach was used in \cite{oliker_advar:05}, and in \cite{glimm-oliker-refl:03}, \cite{glimm-oliker-refl:04}, \cite{galo:04}
it was applied to study problems in geometrical optics.

Concerning regularity of the solution to the Aleksandrov problem, we note that in analytic form, when the measure $\mu$ has a density $m(x)> 0,~x \in \Sn$,
the problem 
reduces to finding the radial function $\p$
satisfying the following equation of Monge-Amp\`{e}re type \cite{Oliker:83}
\begin{equation} \label{MA}
(\p^2 + |\n \p|^2)^{-\frac{n+1}{2}}\p^{1-n}\frac{\det(-\p\n_{ij}\p 
+2\p_i\p_j+\p^2e_{ij})}{\det(e_{ij})}= m ~\mbox{on}~\Sn.
\end{equation}
Here, the gradient $\n$ is computed in the standard metric $e = \sum_{i,j=1}^ne_{ij}du^idu^j$ on $\Sn$, $u^1,...,u^n$ are some local coordinates on $\Sn$, $\n_{ij}\p, ~i,j=1,...,n,$
are the second covariant derivatives in metric $e$, while $\p_i=\frac{\partial \p}{\partial u^i}, i=1,...,n$.
Existence of smooth solutions to the  Aleksandrov problem under
appropriate smoothness assumptions
was established for $n=2$ by Pogorelov \cite{P1}, Chapter VIII,  and
for arbitrary $n$ by the author \cite{Oliker:83}. For related results concerning
$C^0$ apriori bounds and stability of solutions see Kagan \cite{Ka:1974}
and Treibergs \cite{tr:1990}. An analogue of Aleksandrov's problem for nonparametric convex hypersurfaces defined over $\Rnn$ was studied
by Aleksandrov in \cite{Al:42}, Pogorelov \cite{P, P1, P3},
Bakelman \cite{Bak} and, more recently, by McCann \cite{McCann:95} and the author \cite{oliker_advar:05}.


The rest of the paper is organized as follows. 
In section \ref{var1} we study a Legendre-like transform between radial and support functions of hypersurfaces in $\calF$, introduce the class of admissible functions and formulate the minimization problem for the
functional $\calQ$. In section \ref{var2} we specialize this minimization problem
to the class of polytopes and prove existence and uniqueness (up to rescaling) of minimizers. In section \ref{gen} we establish existence and uniqueness of a minimizer 
to the minimization problem for $\calQ$ in the class of general closed convex hypersurfaces. Finally, in section \ref{opt0} we describe the
connections between the variational solution of Aleksandrov's problem
and the mass transport theory.

The author thanks the reviewer for careful reading of the manuscript and many helpful comments.
\section{A variational formulation of Aleksandrov's problem}\label{var1}
The key step leading to a variational statement of the problem
is based on the simple but crucial fact that two representations are available for any
$F \in \calF$ and a transformation generalizing the classical Legendre transform connects these two representations. 
Let $F \in \calF$, $\p$ its radial function and 
$r(x) = \p(x)x,~x \in \Sn,$ its position vector.
An alternative representation of $F$ can be given in terms of its
support function defined
as 
\begin{equation}\label{sup h}
h(N) = \sup_{x\in \Sn}\p(x)\la x, N\ra, ~N\in \Sn.
\end{equation}
Since $\calO \in int B(F)$, where $B(F)$ is the convex body bounded by $F$, functions $\p$ and $h$ are both positive on $\Sn$.
It follows from (\ref{sup h}) that
\begin{equation}\label{sup p}
\frac{1}{\p(x)} = \sup_{N\in \Sn}\frac{\la x, N\ra}{h(N)},~x\in \Sn.
\end{equation}
Thus, the functions $\p$ and $h$ provide two descriptions of the same hypersurface $F$ and these two descriptions are related by (\ref{sup h}) and (\ref{sup p}).
We refer to (\ref{sup h}) as the generalized Legendre-Fenchel (LF) transform of $\p$ and 
write $h =LF(\p)$. Similarly, the function $\p$ defined by (\ref{sup p})
is considered as $LF^{-1}(h)$.
\begin{theorem} \label{charact}
Let $\p>0$ and $h> 0$ be two continuous functions on $\Sn$
satisfying (\ref{sup h}) and (\ref{sup p}). Then there exists a unique closed
convex hypersurface $F \in \calF$ with radial function $\p$ and
support function $h$.
\end{theorem}
\begin{proof} Define the function 
$$H(u) = \sup_{x\in \Sn}\p(x)\la x, u\ra, ~u \in \Rn.$$
Obviously, $H(u)=h(u)$ for $u\in \Sn$ and, thus, $H$ is a positively homogeneous extension of order one  of the function $h$ from $\Sn$ to $\Rn$. It is clear that $H$ is also subadditive. 
By a theorem of H. Minkowski it is the
support function of a unique compact convex body in $\Rn$; see \cite{Schneider}, section 1.7. 
The boundary of that body is a closed convex hypersurface $F$. Since $H(u) > 0$ for $u \neq 0$,
it is clear that $F \in \calF$.

Let 
$$
R(v)=\sup\{\lambda \geq 0~|~\lambda v \in B(F)\}~\mbox{for}~v \in \Rn \setminus \{\calO\}
$$
be the radial function of the convex body $B(F)$ (see
\cite{Schneider}, section 1.7.)
Then $R$ is positively homogeneous of order $-1$ and it is defined by its
values on $\Sn$. Furthermore, $R(x)x \in F$ for any $x \in \Sn$.

We show now that $R(x)=\p(x)$ for all $x \in \Sn$. Fix some $N \in \Sn$ and consider the set
$$
G_N=\{x \in \Sn~|~R(x)\la x, N\ra = h(N)\}.
$$
Note that $G_N \neq \emptyset$. Since $R, h > 0$ on $\Sn$, $\la x, N\ra > 0$ for all $x \in G_N$. Then
by (\ref{sup h}) $R(x) \geq \p(x)$ for all $x\in G_N$.
Suppose there exist $\ox \in G_N$ such that  $\p(\ox) < R(\ox)$.
Since $R(\ox)\ox \in F$, the point $\p(\ox)\ox$ is an interior point of the convex body $B(F)$. Then $\p(\ox)\la \ox,u\ra < h(u)$ for any $u\in \Sn$ and,
consequently,
\[
\frac{1}{\p(\ox)} > \sup_{u\in \Sn}\frac{\la \ox,u \ra}{h(u)}.
\]
But this contradicts (\ref{sup p}).
Thus, $\p(x) =R(x)$ for all $x \in G_N$. Since $\bigcup_{N\in \Sn}G_N = \Sn$,
we conclude that $R(x)=\p(x)$ for all $x \in \Sn$.
\end{proof}

In contrast to the classical Legendre transform (conjugation) defined by
\[
f^*(u) =\sup\{\langle u, x\rangle -f(x)~|~x \in \Rnn\},~u\in \Rnn,
\]
where $f:\Rnn\rightarrow {\bar{\Reals}}$ is a closed convex function
(see  \cite{Schneider}, 
section 1.6), the $LF$ transform relates $\p$ and
$h$ multiplicatively  which is inconvenient for defining the required  functional in the Aleksandrov problem. It is more convenient to work with the logarithms of functions $\p$, $h$
and $\xN$, which are all positive for $F\in \calF$. Put 
$$\hp :=\log \p,~\hh := \log h,~\hK (x,N) := \log \xN~\mbox{for}~ x, N \in \Sn~ \mbox{and}~ \xN > 0.$$
 Then, it follows
from (\ref{sup p}), (\ref{sup h}) that
for a
$F\in \calF$ the generalized Gauss map can be defined as 
\begin{equation}\label{gmap1}
\a_F(x) = \{N\in \Sn ~|~\xN > 0 ~\mbox{and}~\hh(N) - \hp(x) = 
\hK (x,N)\},~ x \in \Sn.
\end{equation}
Similarly,
\begin{equation}\label{gmap2}
\a^{-1}_F(N) = 
\{x\in \Sn ~|~\xN > 0 ~\mbox{and}~\hh(N) -  \hp(x) = 
\hK(x,N)\},~ N \in \Sn.
\end{equation}

In this notation the set of pairs of admissible functions (\ref{dual1}) is
\begin{eqnarray}\label{adm1}
\calA=\{(h, \p) \in C(\Sn)\times C(\Sn)~| ~h > 0,~\p > 0,\\ \nonumber
~\hh(N)-\hp(x) \geq \hK(x,N),~(x,N) \in \Sn\times \Sn,~\xN > 0\}\nonumber
\end{eqnarray}
and the functional (\ref{dual0}) is
\begin{equation}\label{adm2}
\calQ[h,\p] = \int_{\Sn}\hh(N)d\s(N) - \int_{\Sn}\hp(x)d\mu(x),~(h,\p) \in \calA,
\end{equation}
where $\mu$ is as in Theorem \ref{AD1}.
The {\it variational problem} we now consider is to determine a pair $(\Ht,\pt) \in \calA$
such that
\begin{equation}\label{adm3}
\calQ[\Ht,\pt]= \inf_{\calA} \calQ[h,\p].
\end{equation}
\section{A variational solution of Aleksandrov's problem for polytopes}\label{var2}
Just as in the classical case, it is convenient to investigate 
the variational problem above first in the class of convex polytopes in $\calF$.
Let $X_K:=\{x_1,...,x_K\}$ be a set of points on $\Sn$ not contained in
one closed hemisphere. In this section the set $X_K$ will remain fixed. Let $l_1,...,l_K$ denote the set of rays originating
at $\calO$ and going through the points $x_1,...,x_K$. We denote by $\calFK$
the subset of $\calF$ consisting of the convex polytopes with vertices possible 
only on the rays $l_1,...,l_K$. 
Let $P\in \calFK$. Denote by $r_1,...,r_K$ the points of intersection
of $l_1,...,l_K$ with $P$.   
It is useful to recall that $P$ is the boundary of the convex hull of the
points $r_1,...,r_K$. Thus, for the polytope $P$ its radial function $\p$
is completely defined by the vector $(\p(x_1),...,\p(x_K))$. Of course, $r_i=\p(x_i)x_i, ~i=1,...,K$.

For the rest of the paper, when speaking of a convex polytope in $\calFK$
with radial function $\p$  and if there is no danger of confusion, we put for brevity $\p_i:=\p(x_i),~i = 1,...,K$. 

 Note also that  for a Borel subset $\o \subset \Sn$ the integral Gauss curvature $\s(\a_P(\o))=0$ if $x_i \not \in \o$ for any $i=1,...,K$.
 
We will need a version of Theorem \ref{charact} for polytopes.
\begin{proposition} \label{theorem 2 prime} Let $X_K$ be as above, 
$p_1>0,...,p_K>0$, and 
$h \in C(S^n),~h > 0$. Suppose 
\begin{equation}\label{sup h prime}
h(N) = \sup_{1\leq i \leq K}p_i\la x_i,N\ra,~N \in \Sn,
\end{equation}
\begin{equation}\label{sup p prime}
\frac{1}{p_i} = \sup_{N\in \Sn}\frac{\la x_i,N\ra}{h(N)},~ i=1,...,K.
\end{equation}
Then there exists a unique convex polytope $P \in \calFK$ with support
function $h$ and radial function $\p$ such that $\p_i=p_i,~i=1,...,K$.
\end{proposition}
\begin{proof} Consider the convex hull $P$ of the points 
$p_1x_1,...,p_Kx_K$. 
Because $x_1,...,x_N$ are not contained in any closed hemisphere of
$\Sn$ and $p_i > 0$ for all $i=1,...,K$ the polytope $P \in \calFK$. 
Let $\p_P(x),~x\in \Sn,$ and $h_P(N),~N\in \Sn,$ be, respectively, the radial and support functions of $P$. Then $\p_P$ and $h_P$ satisfy (\ref{sup h})
 and (\ref{sup p}). It follows from (\ref{sup h}) and (\ref{sup h prime}) that 
$h_P(N) \geq h(N)~\forall N\in \Sn$. On the other hand, for any fixed $N\in \Sn$ and
any $x\in \Sn$ such that $\p_P(x)\la x,N\ra = h_P(N)$ we have
\[h_P(N) = \sum\lambda_ip_i\la x_i,N\ra \leq h(N),~\lambda_i \geq 0, \sum\lambda_i =1,\]
where the sums are taken over the vertices of $P$ whose linear combination
gives $\p_P(x) x$. Thus, $h_P(N) = h(N)~\forall N \in \Sn$. The equalities 
$\p_P(x_i) = p_i$ for $i=1,...,K$ are proved using arguments similar to those
at the end of the proof of Theorem \ref{charact}.  
\end{proof}

For a vector $p=(p_1,...,p_K) \in \RK$, we write $p > 0$ when $p_1 > 0,...,p_K >0$. Also, for $p >0$ put $\hat{p}=(\hat{p}_1,...,\hat{p}_K)$, where $\hat{p_i}:=
\log p_i,~i=1,...,K$.  Define now the admissible set as
\begin{eqnarray}\label{admp1}
\calA_K=\{(h,p) \in C(\Sn)\times \RK~|~ h>0,~p > 0,  \\ \nonumber
~\hh(N)-\hat{p}_i \geq \hK(x_i,N),~x_i\in X_K,~N \in \Sn,~\xiN > 0,~i = 1,...,K\}.
\end{eqnarray} 
Put
\begin{equation}\label{admp2}
\calQK[h,p] = \int_{\Sn}\hh(N)d\s(N) - \sum_{i=1}^K\hat{p}_i\mu_i,~(h,p) \in \calA_K,
\end{equation}
where $\mu_1,...,\mu_K$ are positive numbers.
When $P$ is a convex polytope in $\calFK$ with the  radial function $\p$, we write $Q_K[h,\p]$ setting in (\ref{admp2}) $p = (\p_1,...,\p_K)$.
\begin{theorem}\label{polytope1}
Let $\mu_i > 0, ~i=1,...,K$,
\begin{equation}\label{sum mu}
\sum_{i=1}^K \mu_i = \s(\Sn)
\end{equation}
and for any (solid) polyhedral convex cone $V\neq \Rn$ with vertex at $\calO$,
possibly degenerate,
\begin{equation}\label{cone cond}
\sum \mu_i > \s((\Sn\cap V)^*),
\end{equation}
where the sum is taken over all $i$ such that $x_i \not \in \Sn\cap V$.
Then there exists a polytope $\Pt \in \calFK$ with radial function $\pt$
and support function $\Ht$ such that
\begin{equation} \label{min}
\calQK[\Ht,\pt]=\inf_{\calA_K} \calQK[h,p]=\inf_{\calFK}\calQK[h,\p].
\end{equation}
The pair $(\Ht,\pt)$ is unique up to rescaling $(c\Ht,c\pt)$ with any
$c=const > 0$.
\end{theorem}
\begin{proof}
First, we note that it suffices to look for a minimizer of $\calQK$ over
$\calFK$. Indeed, if $(h,p) \in \calA_K$ is such that
$$\hh(N) - \hat{p}_i > \hK(x_i,N)$$
for some $i$  and all $N \in \Sn$ such that $\xiN > 0$
then we can increase $p_i$ slightly to $p^{\prime}_i$
and
decrease $h$ to $h^{\prime}$ for some $N$ 
so that the pair
 $(h^{\prime},p^{\prime})$ is still in $\calA_K$. The functional $\calQK$ will
not increase under such change. Therefore, it suffices to minimize $\calQK$ only
on pairs in $\calA_K$ such that
\begin{eqnarray} 
\hh(N) = \sup_i[\hat{p}_i +\hK(x_i,N)],~N\in \Sn, \label{admp5}\\
-\hat{p}_i = \sup_{N\in \Sn} [-\hh(N) + \hK(x_i,N)], ~i=1,...,K, \label{admp6}
\end{eqnarray}
with the supremum achieved in (\ref{admp5}) for each $N$ at some $i$ and 
in (\ref{admp6}) for each $i$ for
some $N$. By Proposition \ref{theorem 2 prime} each such pair $(h,p)$ defines a unique convex polytope in $\calFK$ with the support function $h$ and
radial function $\p$ such that $\p_i = p_i,~i=1,...,K$. 

It follows from (\ref{sum mu}) that if a pair $(h,p) \in \calA_K$ and $c$ is a positive constant, then 
$$
\calQK[ch,cp] = \calQK[h,p].
$$
Hence, the search for a minimizer in $\calFK$ can be restricted
to convex polytopes in $\calFK$ such that $\p_i\leq 1$ for all 
$i\in \{1,...,K\}$ and $\p_i= 1$ at least for one $i\in \{1,...,K\}$.
To avoid introduction of additional notation, we will continue to denote by
$\calFK$ the set of convex polytopes satisfying this restriction, while
$\bFK$ will denote
the closure of the set of such polytopes with respect to the metric
$$
dist (P^1,P^2) = \left[\sum_{i=1}^K (\p^1_i-\p^2_i)^2\right]^{1/2}.
$$
The set of vectors $\{(\p_1,...,\p_K)\}$ corresponding to polytopes
in $\bFK$ forms a closed bounded set in $\RK$ while the corresponding
support functions form a compact set in $C(\Sn)$ with repect to the
uniform convergence on $\Sn$. Note that the functional $\calQK$ is continuous on $\calFK$.

\begin{remark}\label{remark1} Since the set $X_K$ is not contained
in a closed hemisphere, it is easy to construct a polytope in $\calFK$
on which 
\begin{equation} \label{main0}
|\calQK[h,\p]| < \infty.
\end{equation}
{\em For example, the polytope
inscribed in $\Sn$ with vertices $x_1,...,x_K$ satisfies (\ref{main0}).}
\end{remark} 
Next, we  show that
the $\inf_{\bFK}\calQK[h,\p]$ is
attained on a closed convex polytope for which $\p_i > 0 ~\forall i =1,...,K$. Since $X_K$ is not contained in a closed hemisphere, on this polytope $h(N) > 0 ~\forall N \in \Sn$ and then $\calQK[h,\p] > -\infty$. Consequently, the Remark 
\ref{remark1} and the continuity of $\calQK$ on $\calFK$ imply
the existence statement of the theorem.
 
Suppose the $\inf_{\bFK}\calQK[h,\p]$ is attained on a convex polytope
$\Pt=(\Ht,\pt)$ and $\calO \in \Pt$. Again, because $X_K$ is not contained in a closed hemisphere, at least one vertex of $\Pt$ must be at $\calO$. The radial functions of
the polytopes in $\bFK$ are uniformly bounded and there exists a  sequence of
polytopes $P^{\ep} = (h^{\ep},\p^{\ep}) \in \calFK$ converging to $\Pt$  as $\ep \rightarrow 0$. Note that for any $\ep< 1$ 
the inequality
$\p^{\ep}_i \leq \ep$ cannot be true for all 
$i=1,...,K$ as $\p^{\ep}_i =1$ for some $i$. The same  holds also
for $\pt_1,...,\pt_K$.
Let $\{x_{1},...,x_{s}\}$ be the set
of directions for which $\p^{\ep}_i \rightarrow 0$ and $\{x_{s+1},...,x_{K}\}$ the set of directions for which  $\p^{\ep}_i \not \rightarrow 0$. We have
\begin{eqnarray}\label{main1}
\calQK [h^{\ep},\p^{\ep}]=\int_{\Sn}\hh^{\ep}(N)d\s(N)- \sum_{i=1}^K\hp_i^{\ep}\mu_i
=\sum_{i=1}^K \left \{\int_{\a_{P^{\ep}}(x_i)}[\hp_i^{\ep} + 
\hK(x_i,N)]d\s(N) -\hp_i^{\ep}\mu_i \right \} \\
=\sum_{i=1}^s\hp_i^{\ep} \left [ \s(\a_{P^{\ep}}(x_i)) -
\mu_i \right ] + \sum_{i=s+1}^K
\hp_i^{\ep} \left [ \s(\a_{P^{\ep}}(x_i))-\mu_i \right ]
+\sum_{i=1}^K\int_{\a_{P^{\ep}}(x_i)}\hK(x_i,N)d\s(N). \label{main2}
\end{eqnarray}
We consider now separately each of the terms in (\ref{main2}), beginning with the first term on the left.  

Any polytope $P^{\ep}$ is the boundary of the convex hull of vertices of $P^{\ep}$ on the rays $l_{s+1},...,l_K$ and of the polytope 
$P_1^{\ep}$ which is the boundary of
the  convex hull of $\calO$ and 
vertices of $P^{\ep}$ on the rays $l_1,...,l_s$. As $\ep \rightarrow 0$ the polytopes $P_1^{\ep}$ contract to $\calO$
while $P^{\ep}$ converges to $\Pt$. Denote by $V$ the polyhedral angle formed by the faces of $\Pt$ adjacent to $\calO$ and observe that 
$V$
is the boundary of the convex hull 
of the rays $l_{s+1},...,l_K$. 
In addition,
\[\sum_{i=1}^s\s(\a_{P^{\ep}}(x_i)) \rightarrow \s(\a_{V}(\calO))=
\s((\Sn\cap V)^*).\]
It follows from (\ref{cone cond})  that for $\ep$ sufficiently small  we must have
\[\sum_{i=1}^s\mu_i > \sum_{i=1}^s\s(\a_{P^{\ep}}(x_i)).\]
On the other hand, for any $i = 1,...,s$ the convex polytopes
with vertices 
$\{r_j=\p_j^{\ep}x_j,~j=i,s+1,...,K\}$
also converge to $\Pt$.
Consequently,
\begin{equation} \label{main4}
\lim_{\ep \rightarrow 0} \sum_{i=1}^s\hp_i^{\ep} \left [ \s(\a_{P^{\ep}}(x_i)) -
\mu_i \right ] =\lim_{\ep \rightarrow 0}( \min_{1\leq i \leq s}\hp_i^{\ep})  \sum_{i=1}^s\left [ \s(\a_{P^{\ep}}(x_i)) -
\mu_i \right ] = +\infty.
\end{equation}

Since $\p^{\ep}_i \geq \d > 0$ for some $\d$ and all $\ep$ and all $i=s+1,...,K$,
for the second sum in (\ref{main2}) we have
\begin{equation}\label{sum 2}
\left | \sum_{i=s+1}^K \hp^{\ep}_i \left [ \s(\a_{P^{\ep}}(x_i))-\mu_i \right ] \right | < \infty ~\mbox{for all}~ \ep.
\end{equation}

Finally, we consider the last term in (\ref{main2}) and estimate its absolute value from above by a bound independent of $\ep$. 
Observe that for each $P^{\ep}$ and any $x_i \in X_K$ the inequality
$\la x_i, N\ra > 0$ holds for all $N\in \a_{P^{\ep}}(x_i)$,
since the origin $\calO$ is strictly inside the convex body bounded by 
$P^{\ep}$.
Using spherical coordinates
on $\Sn$ with the origin at $x_i$, we obtain
\[
\int_{\a_{P^{\ep}}(x_i)}|\hK(x_i,N)|d\s(N) = 
\int_{\la x_i, N\ra \geq 0} |\hK(x_i,N)|\chi_{\a_{P^{\ep}}(x_i)}(N)d\s(N)
\leq \left \{ \begin{array}{ll} 
\mbox{$\nu_{n-1}$} ~\mbox{when}~ n \geq 2\\
\pi/2 ~\mbox{when}~ n =1,
\end{array}\right.
\]
where $\chi_{\a_{P^{\ep}}(x_i)}$ is the characteristic function of the set
$\a_{P^{\ep}}(x_i)$ and $\nu_{n-1}$ is the $(n-1)-$
dimensional volume of $\Snn$. 
Then
\begin{equation}\label{sum 3}
\left |\sum_{i=1}^K\int_{\a_{P^{\ep}}(x_i)}\hK(x_i,N)d\s(N)\right | \leq
\left \{ \begin{array}{ll} 
\mbox{$K\nu_{n-1}$} ~\mbox{when}~ n \geq 2\\
K\pi/2 ~\mbox{when}~ n =1.
\end{array} \right.
\end{equation}

It follows from (\ref{main1})-(\ref{sum 3}) and  remark \ref{remark1} that $\inf_{\bFK}\calQK[h,\p]$ is attained on a polytope
in $\calFK$. This completes the proof of the existence statement in this theorem. We now prove uniqueness 
of the minimizer of $\calQK$ (up to a homothety
with respect to $\calO$).

Let $P^1=(h^1,\p^1)$ and $P^2=(h^2,\p^2)$ be two minimizers of $\calQK$ in
$\calFK$. After rescaling, if needed, we may assume
that  $\p_i^2 \leq \p_i^1,~\forall i=1,...,K$ with the equality holding at
least for one $i$ and strict inequality for some other $i$. 
Then
$$\calQK[h^2,\p^1] < \calQK[h^2,\p^2]$$ and we arrived at a contradiction.
\end{proof}
\begin{theorem}\label{thm5}
 Let $\Pt$ be a minimizer of $\calQK$ in Theorem \ref{polytope1}. Then 
\begin{equation}\label{polytope 2}
\s(\a_{\Pt}(x_i)) = \mu_i ~\forall i = 1,...,K.
\end{equation}
\end{theorem}
\begin{proof} Suppose that (\ref{polytope 2}) is not true
for some $m$. We show 
that in this case the value of the functional $\calQK$ can be reduced.

Note first that because of (\ref{sum mu}) and because 
\[\sum_{i=1}^K \s(\a_{P}(x_i))= \s(\Sn)~\forall P\in \calFK\]
it may be assumed that $\s(\a_{\Pt}(x_m)) > \mu_m$. 

Consider a convex polytope 
$\Pd$ which is the boundary of the convex hull of points
$r_i=\pd_ix_i,~1\leq i\leq K$, where $\pd_i =\pt_i$ when $i\neq m$ and
$\pd_m=\pt_m-\d$, and
$\d>0$ is sufficiently small so that $r_i \in \Pd ~\forall i$ and $\Pd \in \calFK$. This is possible since $\s(\a_{\Pt}(x_m)) > \mu_m >0$.
(Here the restriction $\p_i \leq 1 ~\forall x_i \in X_K$ with the equality achieved at least for one $i$ imposed on polytopes in $\calFK$
in
the proof of Theorem \ref{polytope1} is not imposed.) The support function of $\Pd$ is denoted by $\hd$. 

Clearly, 
\[
\bigcup_{i=1}^K \a_{\Pt}(x_i) = \bigcup_{i=1}^K \a_{\Pd}(x_i)=\Sn
\] 
and $\a_{\Pd}(x_m)\subset \a_{\Pt}(x_m)$, while $\a_{\Pt}(x_i)\subseteq \a_{\Pd}(x_i)$ for $i \neq m$. Put 
$\kappa_{ij}=\a_{\Pd}(x_i)\cap \a_{\Pt}(x_j)$.
Then 
\[\a_{\Pd}(x_m)= \a_{\Pt}(x_m)\setminus \left (\bigcup_{i\neq m}\kappa_{im}\right ),~\a_{\Pd}(x_i)=\a_{\Pt}(x_i)\cup \kappa_{im}~\mbox{for}~
i\neq m.\]
For $i\neq m$ and $N\in \a_{\Pt}(x_i)$ we have
\begin{equation*}
\hhd(N) = \sup_{1\leq j\leq K}[\hpd_j+\hK(x_j,N)]= \hpd_i+\hK(x_i,N)=\sup_{1\leq j\leq K}[\hpt_j+\hK(x_j,N)]=\hHt(N).
\end{equation*}
For $N\in \a_{\Pd}(x_m)$ we have
\begin{eqnarray*}
\hhd(N) - \hHt(N)= \sup_{1\leq j\leq K}[\hpd_j+\hK(x_j,N)]-
\sup_{1\leq j\leq K}[\hpt_j+\hK(x_j,N)]\\\nonumber
= \hpd_m+\hK(x_m,N) - \sup_{1\leq j\leq K}[\hpt_j+\hK(x_j,N)]\leq
\hpd_m-\hpt_m= -\frac{\d}{\pt_m} +o(\d).
\end{eqnarray*}
In addition, for suitable $i$ (depending on N)
\begin{eqnarray*}
\hhd(N) - \hHt(N) = \hpd_i+\hK(x_i,N)-\sup_{1\leq j\leq K}[\hpt_j+
\hK(x_j,N)]\leq 0~\forall N \in \kappa_{im},~i\neq m.
\end{eqnarray*}
Then
\begin{eqnarray}\label{eq3}
\int_{\Sn}[\hhd(N)-\hHt(N)]d \s(N)=\int_{\a_{\Pt}(x_m)}[\hhd(N)-\hHt(N)]d \s(N)\\ \nonumber+ \sum_{i\neq m}\int_{\a_{\Pt}(x_i)}[\hhd(N)-\hHt(N)]d \s(N)
=\int_{\a_{\Pd}(x_m)}[\hhd(N)-\hHt(N)]d \s(N)\\\nonumber +\sum_{i\neq m}\int_{\kappa_{im}}[\hhd(N)-\hHt(N)]d \s(N)\leq 
\left (-\frac{\d}{\pt_m}+o(\d)\right )\s(\a_{\Pd}(x_m)). 
\end{eqnarray}
It is shown in \cite{Ka:1974}
that
\begin{equation*}\label{kag}
\s(\a_{\Pt}(x_m)) =\s(\a_{\Pd}(x_m)) +C\d + o(\d),
\end{equation*}
where the constant $C$ depends only on the polytope $\Pt$.
Using this in  (\ref{eq3}) gives
\begin{equation*}
\int_{\Sn}[\hhd(N)-\hHt(N)]d \s(N) \leq 
-\frac{\d}{\pt_m}\s( \a_{\Pt}(x_m)) + o(\d).
\end{equation*}
Then for sufficiently small $\d > 0$
\begin{equation*}
\calQK[\hd,\pd]-\calQK[\Ht,\pt]\leq -\frac{\d}{\pt_m}\left [\s( \a_{\Pt}(x_m))-\mu_m \right ] +o(\d) < 0.
\end{equation*}
However, this is impossible since $\inf_{\calFK}\calQK = \calQK[\Ht,\pt]$.
\end{proof}
\section{The general case}\label{gen}
We use now the above results for convex polytopes to solve the
variational problem (\ref{adm3}).
\begin{theorem} \label{thmgen} Let $\mu$ be a measure on $\Sn$ satisfying conditions
{\em (i)-(iii)} in Theorem \ref{AD1}. Then there exists a closed convex hypersurface
$\Ft\in \calF$ with support function $\Ht$ and radial function $\pt$ such
that 
\begin{equation}\label{g1}
Q[\Ht,\pt] = \inf_{\calA}Q[\Ht,\pt].
\end{equation}
In addition, for any Borel set $\o \subset \Sn$ 
\begin{equation}\label{gAD}
\s(\a_{\Ft}(\o))=\mu(\o),
\end{equation}
that is, the minimizer $\Ft$ is a solution of the Aleksandrov problem.
\end{theorem}
\begin{proof} Partition $\Sn$ so that $\Sn = \bigcup_{i=1}^K V_i$,
where $V_i$ are Borel subsets of $\Sn$ such that $V_i \cap V_j = \emptyset$ when $i\neq j$ and $diam(V_i) < \d$. Pick a point $x_i \in V_i$,
and put $\mu_i =\mu(V_i)$. It is shown by Aleksandrov in \cite{Al_invariance} that for $\d$ sufficiently small the set $X_K=\{x_1,...,x_K\}$ and the numbers $\mu_1,...,\mu_k$ (or a subset of $X_K$ corresponding to positive $\mu_i$) satisfy
the conditions of Theorem \ref{polytope1}. Therefore, for each sufficiently large $K$ 
there exists a polytope $\Pt^K \in \calFK$ with support function $\Ht^K$
and radial function $\pt^K$ such that
$\calQK[\Ht^K,\pt^K] \leq \calQK[h^K,\p^K]
~\forall (h^K,\p^K)\in \calFK$.
By rescaling, if needed, it can be assumed that the diameters of all
$\Pt^K$ are equal $1$. Therefore, for $\d \rightarrow 0$ the corresponding sequence $\{\Pt^K\}$ contains
a converging subsequence, which we denote again by $\{\Pt^K\}$. The limit
of this subsequence we denote by $\Ft$ and its support and radial functions
respectively by $\Ht$ and $\pt$.

It is also shown in \cite{Al_invariance} that the  origin $\calO$ is strictly inside the
convex body bounded by $\Pt$ and as $\d \rightarrow 0$ the integral Gauss curvatures of
$\Pt^K$ converge weakly to the integral Gauss curvature of $\Ft$. 
This implies (\ref{gAD}).
The above also implies that
$$\calQK [\Ht^K,\pt^K] \rightarrow \calQ [\Ht,\pt].$$

It is known (\cite{Schneider}, section 2.2) that for any $F \in \calF$
\[\s\{N\in \Sn~|~\a_F^{-1}(N)~\mbox{is multivalued}\} =0.\]
It follows from (\ref{polytope 2}) that, in particular, 
\begin{equation}
\sum_{i=1}^K\int_{\a_{\Pt^K(x_i)}}f(\a^{-1}_{\Pt^K}(N)) d \s(N) = \sum_{i=1}^Kf(x)d\mu^K(x),~\forall f \in C(\Sn),
\end{equation}
where 
$$\mu^K(\o) = \sum_{x_i\in \o}\mu_i$$
for any  Borel subset $\o \subset \Sn$. 
By weak continuity of the integral Gauss curvature,
\begin{equation}\label{g2}
\int_{\Sn}f(\a_{\Ft}^{-1}(N)) d \s(N)  = \int_{\Sn}f(x)d\mu(x),
\end{equation}
which implies (\ref{gAD}).

On the other hand, for any pair $(h,\p) \in \calA$
we have
\begin{equation}\label{g3}
\hh(N) - \hp(\a^{-1}_{\Ft}(N)) \geq \hK(\a^{-1}_{\Ft}(N),N) =\hHt(N) -
\hpt(\a^{-1}_{\Ft}(N))~\forall N \in \Sn.
\end{equation}
Integrating the left hand side of this inequality against the measure $d\s(N)$ and using the change
of variable formula (\ref{g2}) we obtain
\[
\int_{\Sn}\hh(N)d \s(N) - \int_{\Sn}\hp(\a^{-1}_{\Ft}(N))d\s(N) =
\int_{\Sn}\hh(N)d\s(N) - \int_{\Sn}\hp(x)d \mu(x) = \calQ[h,\p].
\]
Integrating the right hand side of (\ref{g3}) gives
\[
\int_{\Sn}\hHt(N)d \s(N) - \int_{\Sn}\hpt(\a^{-1}_{\Ft}(N))d\s(N) =
\int_{\Sn}\hHt(N)d\s(N) - \int_{\Sn}\hpt(x)d \mu(x) = \calQ[\Ht,\pt].
\]
Thus $\calQ[\Ht,\pt] \leq \calQ[h,\p] ~\forall (h,\p) \in \calA$. This proves (\ref{g1}).
\end{proof}
\begin{theorem} \label{uniqueness} The minimizing pair in Theorem \ref{thmgen} is unique
up to rescaling with any constant $c > 0$.
\end{theorem}
\begin{proof} Let $(h_i,\p_i),~i=1,2$, be two minimizing pairs which
are not constant multiples of each other. Rescale, if needed, one
of the pairs so that
\begin{equation*}
\p_2(x) \leq \p_1(x) ~\forall x \in \Sn
\end{equation*}
with the equality holding for at least one $x$.
It follows from (\ref{sup h}) that  
$h_2(N) \leq h_1(N) ~\forall N \in \Sn$. Assuming that $\p_2 \not \equiv \p_1$, it is clear that there exists $\bar{N}$ such that $h_2(\bar{N}) < h_1(\bar{N})$. Since the support function of
a convex body is continuous, this inequality holds for some set of positive measure on $\Sn$. Then
\[Q[h_2,\p_1] < Q[h_1,\p_1].\]
This contradicts the minimality of $(h_1,\p_1)$. Therefore,
$h_2(N) = h_1(N)$ almost everywhere on $\Sn$. 

Let us show that $h_2(N) = h_1(N)~\forall N \in \Sn$. Extend both functions from $\Sn$ to the entire $\Rn$ as positively homogeneous functions of order 1 and denote these extensions by $H_1$ and $H_2$.
The functions $H_1$ and $H_2$ are locally Lipschitz and differentiable almost everywhere \cite{Schneider}, section 1.5. Therefore, $d H_2 -d H_1 =0$ almost everywhere,
that is, $H_2-H_1 =const$ almost everywhere. Then $H_2-H_1=0$
almost everywhere. Since $H_2-H_1$ is locally Lipschitz, we conclude that $H_2-H_1=0$ everywhere. This implies
that $h_2-h_1 \equiv 0$ on $\Sn$ and consequently $\p_1\equiv\p_2$.
\end{proof}
\begin{remark} The convergence of minimizers $\Pt^K$ to a minimizer of
$\calQ[h,\p]$ combined with the uniqueness of the minimizer provides a constructive way for finding the solution to Aleksandrov's problem numerically.
\end{remark}
\section{Aleksandrov's problem and optimal transport on $\Sn$}\label{opt0}
\subsection{Connection with the problem of Monge}\label{opt0-1}
Consider the problem of finding a map $\th_0 \in \Th$  such that
\begin{equation}\label{monge}
\sup_{\th \in \Th} \int_{\Sn} c(\th^{-1}(N),N)d\s(N) =  \int_{\Sn} c(\th_0^
{-1}(N),N)d\s(N),
\end{equation}
where $\Th$ is defined by (\ref{theta1})-(\ref{g4})  and $c$  by (\ref{cost}) in the introduction. Note that $\Th \neq \emptyset$. For example,
the generalized Gauss map  $\a_{\Ft}$ of the minimizer in Theorem \ref{thmgen} is in 
$\Theta$. The required properties of $\a_{\Ft}$ follow from known results in convexity
theory (see, for example, \cite{Schneider}, section 2.2) and the change of variable formula (\ref{g2}). 

The maximization  problem (\ref{monge}) can be viewed as a variant of the celebrated Monge problem \cite{Monge81}. The original problem of Monge is formulated for two measurable sets $U$ and $V$ 
in Euclidean space $\Rnn$ and two Borel measures $m_U$ and $m_V$ with equal total masses defined on $U$ and $V$, respectively. It consists in finding
a map $T_0$ among  all measurable maps $\{T\}$ of $U$ onto $V$ such that
$m_U(T^{-1}(B)) = m_V(B)$ for all Borel sets $B \subseteq V$ (that is,
$T$ pushes $m_U$ forward to $m_V$) , and 
\[\inf_{\{T\}}\int_Ul(x,T(x))dm_U = \int_Ul(x,T_0(x))dm_U, 
\]
where $l$ is the Euclidean distance; see, for example, \cite{KA}, Ch. VIII,
section 4, or \cite{vilani:04}. 
In our case the problem (\ref{monge}) is
considered on $\Sn$ and the cost function $c(x,N)$  is nonlinear and allowed
to assume infinite values.
\begin{theorem}\label{monge1} The problem  (\ref{monge}) admits a solution $\th_0$ and any such solution satisfies
$\th_0^{-1}= \a_{\Ft}^{-1}$ almost everywhere on $\Sn$ with respect to $\s~ (\s-\mbox{a.e.})$. Here, $\Ft\in \calF$  is the unique (up to rescaling) convex hypersurface defined by the unique (up to rescaling) minimizing pair $(\pt,\Ht)$
of the functional $Q$ in Theorem \ref{thmgen}. Furthermore, 
\begin{equation}\label{dual0-0}
\calQ[\Ht,\pt]=\int_{\Sn} c(\a_{\Ft}^{-1}(N),N)d\s(N)=\sup_{\th \in \Th} \int_{\Sn} c(\th^{-1}(N),N)d\s(N).
\end{equation}
In addition, for any Borel set $V \subset \Sn$ for which $\mu(V) > 0$ we have 
$\th_0(V)=\a_{\Ft}(V)$. 
\end{theorem}
\begin{proof} Let, as before, $(\Ht,\pt)$ be the pair in $\calA$ defining
$\Ft$. For an arbitrary pair $(h,\p) \in \calA$ and any $\th \in \Th$ we have for almost all $N\in\Sn~\mbox{such that}~\langle \th^{-1}(N),N \rangle > 0$
\begin{equation}\label{monge2}
\hh(N) - \hp(\th^{-1}(N)) \geq \hK(\th^{-1}(N),N).
\end{equation}
Integrating against $d \s(N)$ and using 
(\ref{g4}), we obtain
\begin{equation}\label{monge3}
Q[h,\p(\th^{-1})] \geq \int_{\{N\in\Sn~|~\langle \th^{-1}(N),N \rangle > 0\}}\hK(\th^{-1}(N),N)d\s(N) \geq \int_{\Sn}c(\th^{-1}(N),N)d\s(N).
\end{equation}
On the other hand, if $h=\Ht,~\p =\pt$ and $\th=\a_{\Ft}$ then 
\begin{equation}\label{monge4}
\hHt(N) - \hpt(\a_{\Ft}^{-1}(N)) = \hK(\a_{\Ft}^{-1}(N),N)~\forall N\in \Sn.
\end{equation}
Since $\a_{\Ft} \in \Th$, the right hand side in (\ref{monge3}) attains its supremum on $\th_0$ such that 
\begin{equation}\label{monge5}
\langle \a_{\Ft}^{-1}(N),N \rangle = \langle \th_0^{-1}(N),N\rangle
~\mbox{$\s-$a.e. on} ~\Sn.
\end{equation}

Let us show that $\th_0^{-1}=\a_{\Ft}^{-1} ~\mbox{$\s-$a.e. on}~ \Sn$. It follows from (\ref{monge4}) and (\ref{monge5}) that
\begin{equation}\label{monge6}
\hpt(\th_0^{-1}(N)) = \hpt(\a_{\Ft}^{-1}(N))~\mbox{$\s-$a.e. on} ~\Sn. 
\end{equation}
Let $M$ be the union of the sets where either $\a_{\Ft}^{-1}$ or
$\th_0^{-1}$ is multivalued and let $N \in \Sn \setminus M$. 
By Theorem \ref{uniqueness} the map $\a_{\Ft}$ is defined uniquely and then by (\ref{monge5}) and (\ref{monge6})
\[\langle \pt (\th_0^{-1}(N))\th_0^{-1}(N),N\rangle =
\langle \pt (\a_{\Ft}^{-1}(N))\a_{\Ft}^{-1}(N),N\rangle.
\]
This means that the points $\pt (\th_0^{-1}(N))\th_0^{-1}(N)$ and  $\pt (\a_{\Ft}^{-1}(N))\a_{\Ft}^{-1}(N)$ lie
on the same hyperplane with normal $N$ supporting to $\Ft$. Since $\Ft$ is convex, it contains the linear segment joining these two points and
this linear segment is also contained in the same hyperplane. If
$\a_{\Ft}^{-1}(N)\neq \th_0^{-1}(N)$ then this segment does not reduce to a point and then $\a_{\Ft}^{-1}(N)$ is not single valued.
This contradicts the choice of $N$.
Therefore, $\th_0^{-1}=\a_{\Ft}^{-1}~\mbox{$\s-$a.e. on} ~\Sn$.
This and (\ref{monge4}) imply that we have equalities
in (\ref{monge3}) when $\th^{-1}=\a_{\Ft}^{-1}$. Consequently, (\ref{dual0})
holds. 

Now we prove the last statement of the theorem.
Let $V\subset \Sn$ be a Borel set and $\mu(V) > 0$. By (\ref{gAD}) 
$\s(\a_{\Ft}(V))=\mu(V) > 0$ and then
$\a_{\Ft}^{-1}(\a_{\Ft}(V))=\th_0^{-1}(\a_{\Ft}(V))$ $\s-$ a.e. in 
$\a_{\Ft}(V)$. Therefore,
$$V= V\cap\a_{\Ft}^{-1}(\a_{\Ft}(V))=V\cap\th_0^{-1}(\a_{\Ft}(V)).$$
Applying $\th_0$ to both sides of the last equality, we obtain $\th_0(V) = \a_{\Ft}(V).$
\end{proof}
\begin{remark} A curious feature of the problem (\ref{monge}) is
that it is, in fact, an intrinsic problem on the sphere $\Sn$ in the
sense that the class of admissible maps $\Th$ and the cost
functional are expressed in terms of the logarithm of the cosine
of the geodesic distance on $\Sn$. Though the condition (\ref{eq1}) in Theorem \ref{AD1} is formulated here for convenience as an extrinsic condition it is in fact intrinsic as the notion of polarity(=duality) in $\Sn$ and (\ref{eq1}) can be stated intrinsically \cite{tr:1990}. Remarkably, the optimal
map is (almost everywhere) the generalized Gauss map of a closed convex
hypersurface in $\Rn$. Thus, the solution of the optimization problem (\ref{monge}) on $\Sn$ is a
solution of the problem of embedding $\Sn$ into $\Rn$ with prescribed integral Gauss curvature.
\end{remark}
\subsection{A representation of the generalized Gauss map}\label{gmap}
For the optimal pair $(\pt,\Ht)$ and the corresponding convex hypersurface $\Ft \in \calF$ the equation (\ref{gmap1}) defines implicitly the generalized Gauss map $\a_{\Ft}$. In this section we give an explict representation for this map in a form resembling the
representation of optimal maps between subsets in Euclidean space with convex or concave cost functions \cite{Brenier}, \cite{Brenier1}, \cite{Gangbo/McCann:95}, \cite{Gangbo/McCann:96}, \cite{Caf_alloc:96}. However, in our case we give a representation for
$\a_{\Ft}$ valid at all points of $\Sn$ and the $\mbox{spt} \mu$ is allowed
 to be atomic. This stands in contrast with usual assumptions under
which representations of optimal maps are derived \cite{Gangbo/McCann:96}. 
In order to state our result we need some definitions.
\begin{definition} Let $f \in C(S^n)$, $x_0 \in \Sn$ and 
$T \Sn_{x_0}$ the tangent space to $\Sn$ at $x_0$.  The set
\[
\partial f(x_0) = \{ v \in  T \Sn_{x_0}~|~ f(x)\langle -v +f(x_0)x_0,x\rangle \leq f^2(x_0)~\forall x \in \Sn\}
\]
is called the subdifferential of $f$ at $x_0$.
\end{definition}
The following proposition clarifies the geometric meaning of a subdifferential in the case of the radial function $\p$ of a 
convex hypersurface $F\in \calF$.
\begin{proposition}\label{gmap3}
Let $F\in \calF$ and $\p$ its radial function. For an arbitrary $x\in \Sn$ denote by $M(x)$ an outward normal to a hyperplane supporting to $F$ at $r(x)=\p(x)x$ rescaled (if needed) so that $\langle M(x),x\rangle = \p(x)$. Then
$v=-M(x) + \p(x)x \in \partial \p(x)$. Conversely, for each $v \in \partial \p(x)$ there exists a unique hyperplane supporting to $F$ at $r(x)$ with outward normal $M(x)$ such that $v = -M(x)+\p(x)x$.
\end{proposition}
\begin{proof} Fix some $x_0 \in \Sn$ and let $P$ be a hyperplane supporting to $F$ at $r_0 = \p(x_0)x_0$. The set of supporting hyperplanes to $F$ is not empty at any point of $F$ 
and such $P$ exists. Let $M_0$ be the outward normal to $P$. Denote by $M_0^{\perp}$ the projection of $M_0$ onto the one-dimensional subspace generated by $x_0$ and let $M_0^{\top} = M_0 -M_0^{\perp}$.
Note that $\langle M_0,x_0 \rangle > 0$ because the origin $\calO$ is strictly inside the convex body bounded by $F$ and $M_0$ is an outward normal. 
Put 
$$M_0^{\prime} = \frac{\p(x_0)M_0}{\langle M_0,x_0\rangle}.$$ Then 
$\langle M_0^{\prime}, x_0\rangle = \p(x_0)$. Since $P$ is supporting to $F$ at $r_0$
we have
\[
\langle\p(x)x, M_0^{\prime} \rangle = 
\p(x)\langle x, M_0^{\prime\top}+M_0^{\prime\perp}\rangle \leq \langle\p(x_0)x_0, M_0^{\prime} \rangle =\p^2(x_0)~\forall x\in \Sn.
\]
Thus,
$v := -M_0^{\prime\top}\in \partial \p(x_0)$. This also shows that  for any $x\in \Sn$ the projection on $T\Sn_{x}$ of an appropriately rescaled outward normal to any supporting hyperplane to $F$  at $\p(x)x$ gives a uniquely defined  element in $\partial \p(x)$. 

Conversely, let $v \in \partial \p_0$. Consider a hyperplane with normal vector $M_0 = -v + \p(x_0)x_0$ containing the point $r_0$.
It follows immediately from the definition of $v$ that this plane is supporting to $F$ at $r_0$.
\end{proof}
\begin{remark}
At points of differentiability of $F$ where there exists only one supporting hyperplane the function $\p$ is differentiable and
$\partial \p(x) = \{\mbox{\/{\em grad}}\p(x)\}$ with $\mbox{\/{\em grad}} \p(x)\in T\Sn_x$, where the gradient is computed with respect to the standard metric of $\Sn$. 
By Rademacher's theorem almost all points on $F$ (in the sense of Lebesgue measure
on $F$) and consequently on $\Sn$ (since it is a radial projection of $F$) are points of differentiability \cite{Schneider}, notes to section 1.5.
\end{remark}
\begin{corollary}
Let $F\in \calF$ and $\p$ its radial function. Then 
for each $x\in \Sn$ and each $N(x) \in \a_F(x)$ there exists a unique
$v(x) \in \partial \p(x)$ such that
\begin{equation}\label{gmap4}
N(x) = \frac{-v(x)+\p(x)x}{\sqrt{v^2(x)+\p^2(x)}}
\end{equation}
and
\begin{equation}\label{gmap5}
\a_F(x) = \left \{\frac{-v(x)+\p(x)x}{\sqrt{v^2(x)+\p^2(x)}} ,~v(x) \in \partial \p(x)\right\}.
\end{equation}
\end{corollary}
\subsection{Connection with the problem of Kantorovich}
In this section we show that in the language of the optimal mass transport theory the minimization problem (\ref{adm3}) can be viewed as the dual of the following primal problem. Denote by $\Gam$ the set of joint Borel measures
on $\Sn\times\Sn$ with marginals $\mu$ and $\s$, where $\mu$ is as in Theorem \ref{AD1} and $\s$, as before, the standard Lebesgue measure on $\Sn$. Thus, any $\g \in \Gam$ satisfies 
\begin{equation}\label{primal1}
\g[U,\Sn] = \mu[U]~\mbox{and}~\g[\Sn,U]=\s[U]~\mbox{for any Borel}~U \subset \Sn.
\end{equation}
The {\it primal} problem in this setting is to determine
a $\tg \in \Gam$ such that
\begin{equation}\label{primal2}
\calC[\tg] = \sup_{\Gam}\calC[\g]
\end{equation}
where $\calC[\g]$ was defined by (\ref{primal}).
The duality relation between the problems (\ref{primal2}) and  (\ref{adm3}) is given by the following
\begin{theorem}\label{duality1}
Let $\Gam$ and $\calC$ be as above. Then there exists a $\tg$ satisfying
(\ref{primal2}) and
\begin{equation}\label{duality3}
\calQ[\Ht,\pt] =\calC[\tg],
\end{equation}
where $(\Ht,\pt)$ is the optimal solution of problem (\ref{monge1}).
\end{theorem}
\begin{proof} It is clear that the Gauss map of the optimal solution $\Ft=(\Ht,\pt)$ in Theorem \ref{monge1}
gives a measure $\g_o \in \Gam$. In fact, any map from the set $\Theta$ defined in section \ref{opt0-1} gives a measure in $\Gam$. To see this, note that if $\th \in \Th$ then one can take
$\g_{\th}[U,V] = \s [\th(U)\cap V]$ for any Borel sets $U,V$ on $\Sn$. Then by (\ref{g4})
$$\calC[\g_{\th}]=
 \int_{\Sn}c(\th^{-1}(N),N)d\s(N).$$
Consequently, taking into account (\ref{dual0-0}), we get
\begin{equation}\label{duality2}
\sup_{\Gam}\calC[\g] \geq \sup_{\th \in \Th} \int_{\Sn} c(\th^{-1}(N),N)d\s(N)
=\int_{\Sn} c(\a_{\Ft}^{-1}(N),N)d\s(N)=\calQ[\Ht,\pt].
\end{equation}

We prove now the reverse inequality. Since $c(x,N) \leq 0 ~\forall (x,N) \in \Sn \times \Sn$, it is clear
that in the maximization problem (\ref{primal2}) it suffices to consider only $\g \in \Gam$ with support
\begin{equation}\label{duality4}
\spt
 \g \subset  \{(x,N) \in \Sn \times \Sn~ |~\xN \geq 0 \}. 
\end{equation}
Denote by $\Gamp$ the subset of measures in $\Gam$ satisfying (\ref{duality4}). 
For any pair $(\hh,\hp)$ from the set $\calA$ of admissible functions (\ref{adm1}) and any $\g \in \Gamp$ we obtain, taking into
account (\ref{adm2}),
\begin{equation}\label{duality5}
\calC[\g]\leq \int_{\Sn}\int_{\Sn}[\hh(N)-\hp(x)]\g(dx,dN) =
\int_{\Sn}\hh(N)\g(\Sn,dN)-\int_{\Sn}\hp(x)\g(dx,\Sn)=\calQ[h,\p].
\end{equation}
It follows from (\ref{duality2}) and (\ref{duality5}) that the
supremum in (\ref{primal2}) is attained on $\tg$ corresponding
to $(\Ht,\pt)$, that is, $\g_o =\tg$ and (\ref{duality3}) holds.
\end{proof}
\bibliographystyle{plain} 
\bibliography{../../../Nov.02}
\end{document}